# On Sensitivity of Time Step for Dynamic Analysis of Bridges under Moving Loads

**Sajad Ahmad Hamidi[1]**

[1] Department of Civil and Environmental Engineering, University of Wisconsin-Milwaukee,
hamidi@uwm.edu

**Abstract:**
In the analysis of structures under dynamic loads, selection of a proper time step has a great influence to reach exact results. In this research determination of proper time step in dynamic analysis of railway bridges under high speed moving loads is considered. Dynamic responses of four simple span steel bridges with 10, 15, 20 and 25 meter length to moving trains with speed from 100 to 400 km/h and axle distances from 13 to 23 meter are calculated considering different time steps in analysis. The results indicate that by increase in moving speed of vehicle (increase in loading speed) the length of proper time step for dynamic analysis is reduced. In contrast by increase in span length (increase in bridge vibration period) longer time steps can be used in dynamic analysis. In this research by investigation of dynamic analysis results, an equation is suggested for determination of proper time step for dynamic analysis of bridges under moving loads.

**Introduction**
Analysis of bridges under moving vehicles based on dynamic behavior of loading is impossible by common static methods and in recent decade's researchers proposed different analytical and mathematical solution to deal with the complicacy of the problem.

The basis of all methods is integration with respect to time and determination of structural responses history during the motion of loads. Integration methods in dynamic of structures are from most efficient methods in analysis of a structure under a dynamic load. In most of integration methods like Newmark and Wilson methods for convergence and improvement the accuracy of results, for time step used in analysis, some constraints must be considered. For example in an analysis with Newmark method for accurate results the time step must be less than the 0.55T where the T is period of first mode of the structure's viberation. In finite difference and Fox-Goodwin methods the proper time steps for convergence of analysis are 0.32T and 0.39T respectively [1].

In weighted residual which is proposed by Razavi and Ghasemiaeh in 2007, the maximum time step for stability of the method is 1.24T [2].

Determination of time step for a dynamic analysis is important not only for stability and convergence of an analysis but also for accuracy of the results. When the amount of error transmitted from one step to the next one is high, after some iteration the analysis may diverge. In 1973 Wilson and Bathe implanted a comprehensive study on stability and accuracy of integration methods in dynamic of structures [3].

On the other hand in some methods of dynamic analysis, the selection of a time step has no visible effect on stability of the calculation but may cause inaccurate results. In this situation as the analysis ensue to the result it is probable the error in results not to be considered. So even in stable methods and analysis it is important to pay attention on the effect of time step on the accuracy of the results. Especially when the dynamic loads imposed on a structure have a high velocity the selection of time step for analysis is from greater importance. If the time step is not sufficiently small, the responses at the beginning and at the end of a time step which are considered as results of analysis may differ dramatically with the responses during the time step. For example if the maximum of a response occurs in the midpoint of a time step, it may be lost during the analysis because the analysis runs only on the beginning and the end of the time step. The time step must be determined so that the error is acceptable.

On the other hand unnecessarily small time step leads to extra time and cost for analysis of complicated problems. So determination of a proper time step is indeed the choice of optimum time step which reduces the error to an acceptable amount and at the same time doesn't increase the number of iterations, cost and time of an analysis unnecessarily.

In common methods for bridge design, in order to account for dynamic effects of vehicle load, traffic load is assumed as a static load increased by implementing the impact factor. The impact factor (I) is defined based on maximum value of dynamic and static responses:

$$I = \frac{D_{dyn} - D_{st}}{D_{st}} \Rightarrow \frac{D_{dyn}}{D_{st}} = 1 + I \qquad (1)$$

As the impact factor based on deflection is greater than those based on other responses like acceleration [4] in this paper dynamic and static deflection at the midpoint of bridges considered the calculation of impact factor. In [5] different relations for Impact Factor based on current codes are available.

In this research for determination of proper time step in dynamic analysis of bridges under moving loads, impact factor based on midpoint deflection is calculated for different situations of bridge and train by considering different time steps in analysis. By comparison of results for different time steps the time step which leads to accurate results in different situation will be determined.

## 3- Dynamic models

In dynamic analysis of bridge under moving vehicles as the responses of bridge like deflection and acceleration in vertical direction is dominant in most cases two dimensional modeling is enough. One axis is in bridge direction which is the direction of load motion too and another is vertical direction which is the direction of loading. The equation of motion for the bridge, showed in Figure 1, is:

$$m\ddot{u} + cIu'''' + EIu'''' = P\sum_{j=1}^{N} \delta(x - v(t - t_j)) * \left[ H(t - t_j) - H(t - t_j - \frac{L}{v}) \right]. \qquad (2)$$

In this equation $u'$ is the derivative of u relative to x coordinate and $\dot{u}$ is derivative of u relative to time. Other symbols in the equation are: m: mass of unit length of beam,

u(x,t): vertical displacement of beam, C: damping coefficient, E: elasticity module I: moment inertia of beam, δ: delta function, H(t): Heavy side function, tj=(j-1)d/v: arriving time of the jth load at the beam and N: the number of moving loads.

By considering initial and boundary conditions and rewriting u(x; t) based on generalized coordinates and mode shape functions Eq. (3) is achieved. Details can be found in [6].

$$\ddot{q}_n(t+\delta t) + 2\xi\omega_n \dot{q}_n(t+\delta t) + \omega^2{}_n q_n(t+\delta t) = F_n(t+\delta t) \tag{3}$$

in which values $\omega_n$ and $F_n(t)$ are calculated as follows:

$$\omega_n = n^2\pi^2 \sqrt{\frac{EI}{mL^4}}$$

$$F_n(t+\delta t) = \frac{2P}{mL}\sum_{j=1}^{N}\left[\sin(\frac{n\pi v(t+\delta t-t_j)}{L})H(t+\delta t-t_j)\right] \tag{4}$$

$$+ \frac{2P}{mL}\sum\left[(-1)^{n+1}\sin(\frac{n\pi v(t+\delta t-t_j-\frac{L}{v})}{L}).H(t+\delta t-t_j-\frac{L}{v})\right]$$

Using Duhamel's integral method the above deferential equation is solved. For the sake of brevity we refer the readers for more detail to [4, 6].

### 3.1. Train modeling

Various models from the simple model such as concentrated loads, or complicated ones such as the sprung mass model, have been used in vehicle modeling. In many studies it has been proved that increasing the modeling details of the train just increases the precision of vehicle response calculations and is not effective in bridges [6, 7]. Therefore when the aim is studying bridge reactions and construction aspects, considering fewer details for modeling suffices. In this research, concentrated loads model has been used to analyze bridge dynamic reactions. The train's specifications have been specified with regard to high speed trains whose samples are presented in References [6]. Each axle load is 20 tons and each train consists of 10 axles with equal distances. In order to study the effect of axle distance on the bridge dynamic reactions, 12 axle distances from 13 to 23 m are considered. Train velocity range is considered from 100 to 400 km/h. The step for velocity increase for each dynamic analysis is 2.5 m/s so 34 different velocities have been analyzed for each axle distance. Table 1 shows the train specifications in dynamic analysis of the research.

Table 1: Train specifications in dynamic analysis

| Train specifications | Variations | Number | Number of dynamic analysis |
|---|---|---|---|
| Load per axle (ton) | 20 | 1 | 408 |

| Number of axles | 10 | 1 | |
|---|---|---|---|
| Axle distance (m) | 13; 14; : : : ; 24 | 12 | |
| Speed (km/h) | 109; 118; : : : ; 406 | 34 | |

## 3.2. Bridge modeling

Most of steel bridges are single span with simple supports; therefore considering a bridge as a simply supported beam (Fig. 4) is the most common method since the main vibrations and displacements occur in a vertical direction. It suffices to consider the problem as two dimensional in order to study dynamic responses.

In this research four simply supported steel bridges with 10, 15, 20 and 25 m span lengths are studied. The dynamic characteristics of these bridges are presented in Fig. 5 and Table 2. For each Bridge, dynamic analysis for different train specifications and different time steps are carried out in this research.

Table 2: Bridge specifications

| Span length (m) | First mode flexural frequency (Hz) | First mode flextural period (s) |
|---|---|---|
| 10 | 12 | 0.083 |
| 15 | 8 | 0.123 |
| 20 | 6 | 0.166 |
| 25 | 4.8 | 0.208 |

## 4- Analysis method

In this research, an analytical method has been used, considering the bridge as a simple beam and train axles as concentrated loads, in order to calculate the bridge responses to train movement.

By these assumptions for modeling and based on Equations 2_4, a program in Matlab environment has been developed that can calculate deflections at various points of the bridge while train axles are passing over, and also calculate the maximum amount of deflection at the midpoint of bridge which is necessary for calculation of impact factor. More details about the program are available in [4].

As mentioned in part 1, in most of integration methods in dynamic of structures the proper time step for convergence of the analysis is a coefficient of the period of first mode of structure. According to this idea at the first step of this research the time step considered 0.1 of the greatest period of flexural mode of above mentioned bridges, which is 0.025 second. Then for investigation of the effects of time step on results of dynamic analysis, each analysis was repeated with considering different time steps: 0.05, 0.025, 0.015, 0.01, 0.005 and 0.0025.

By comparison the results of dynamic analysis for different time step, the proper time step which leads to acceptable accuracy in responses is chosen. If analyses show in some cases 0.05 second results to accurate responses, greater time step must be considered in dynamic analysis to find the greatest time step which leads to acceptable accuracy in results. On the other hand if the results for dynamic analyses with 0.0025 second are not exact in some cases the analysis for those cases must be repeated with smaller time steps.

It is worth mentioned the code developed in MATLAB was validated by comparison with dynamic responses of bridges. More detail about the validation of the model is available in [4].

## 5- Results

As mentioned in previous section time step is considered 0.025 second and midpoint deflection and impact factor are calculated. Some results are shown in Figs 2-5.

In some cases like that in impact factor for bridge with 15 m length under moving train with 13 m axle distance the diagram decreases when speed increases from 360 to 375 kilometer per hour. According to equation 6 the impact factor for speed 375 km/h must be the maximum as a result of coincidence of loading frequency with frequency of bridge vibration. This is shown by arrow in the figure 3.

$$V = 375 km/h = 104.16 m/s \rightarrow f_l = \frac{V}{d} = \frac{104.16}{13} = 8.01 \approx 8 \qquad (7)$$

So the value of impact factor for this situation is not accurate and for determination of exact value for impact factor the dynamic analysis of bridge repeated with considering time step 0.015 second. The results of both analyses are presented in figure 6. it is obvious the reduction of time step in analysis leads to increase in accuracy of results and the maximum occurs at the speed of 375 km/h as expected.

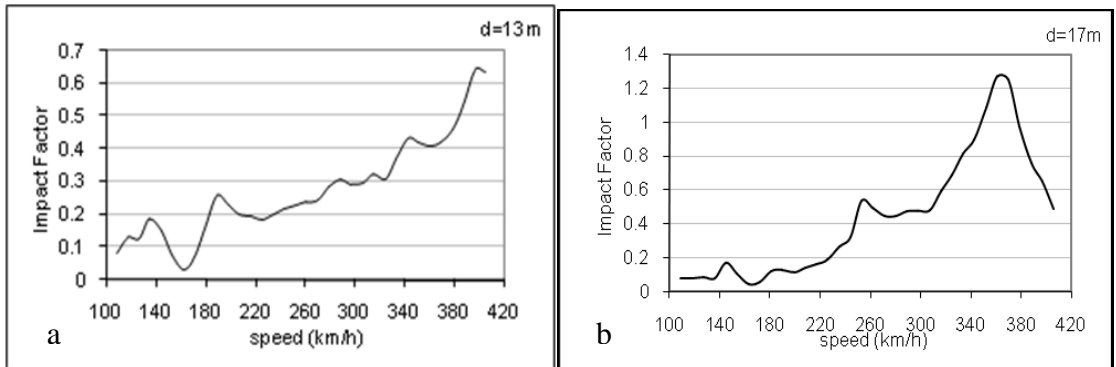

Figure 2: Impact factor for bridge with 10 m span length and train with 13 and 17 m axle distances

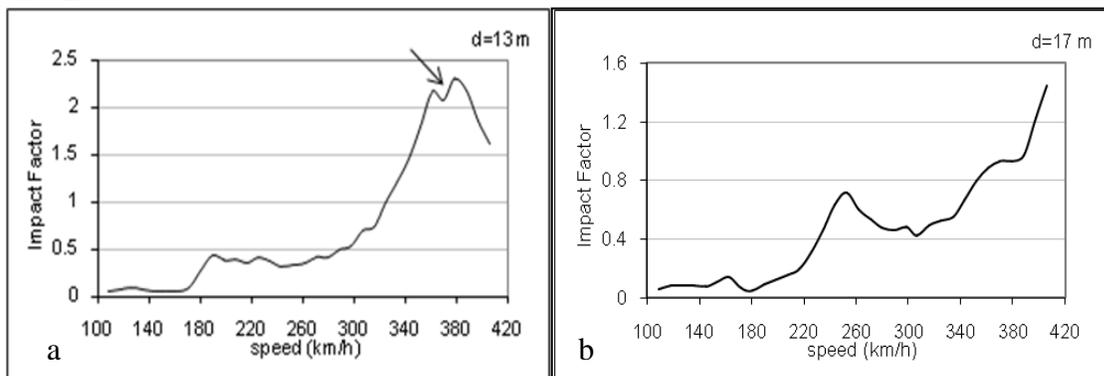

Figure 3: Impact factor for bridge with 15 m span length and train with 13 and 17 m axle distances

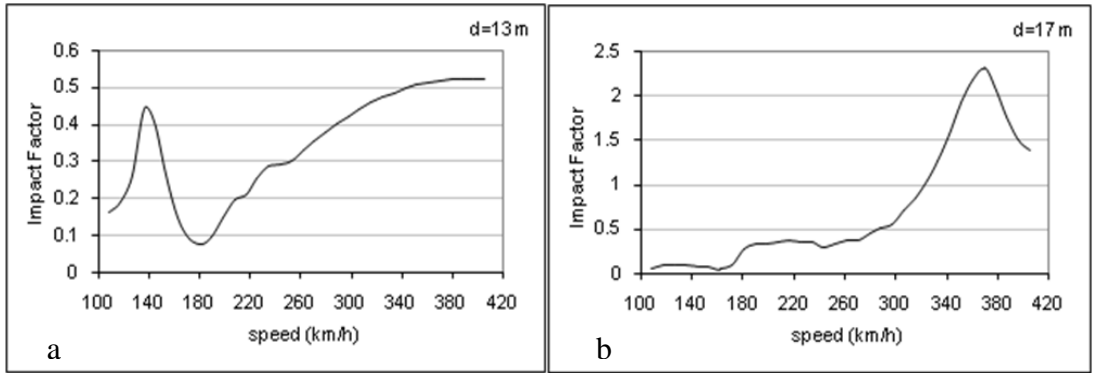
Figure 4: Impact factor for bridge with 20 m span length and train with 13 and 17 m axle distances

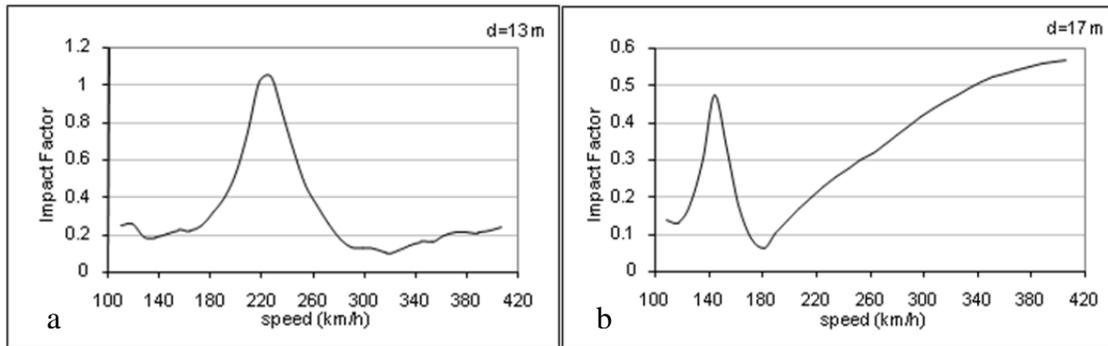
Figure 5: Impact factor for bridge with 25 m span length and train with 13 and 17 m axle distances.

Although decrease in time step leads to more accuracy in results of a dynamic analysis but at the same time, the number of calculations, time and cost will be increased. So determination of time step so that the dynamic analysis results to accurate responses in the minimum possible time and cost for different condition is useful. If a predefined proper time step exists there is no need for checking the accuracy of results nor trail and error to find the proper time step for each analysis. This results in a huge reduction in efforts and consequently the cost of a dynamic analysis in research and engineering.

For determination of proper time step in this research, dynamic analyses are carried out for different time steps, 0.05, 0.025, 0.015, 0.01, 0.005 and 0.0025. By calculation of maximum midpoint deflection and impact factor for different condition the parameters which affect the amount of proper time step in analysis of bridge under moving load are determined.

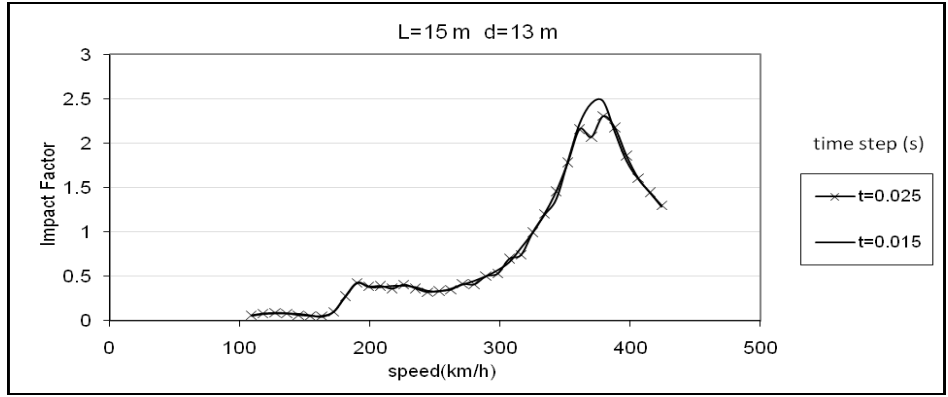

Figure 6: Impact factor for bridge with 15 m span length and train with 13 m axle distance and different time steps

### *5-1- Moving speed*

In figures 7 to 10 impact factors for different speeds and different time steps are shown. By investigation of results, some of them are shown in following figures, it is obvious that in higher speeds, higher difference in impact factor calculated by different time step is more probable.  In higher speed for accurate result dynamic analysis must be implemented with a smaller time steps. Indeed by increase in speed of moving loads, the dynamic responses of structure change rapidly. If the time step is not small enough it is possible a maximum value occurs in the middle of a time step and is not shown in an output of analysis. So it is important to reduce the time step by increase in velocity of the moving loads. In conclusion, there is an inverse proportion between proper time step in dynamic analysis and velocity of moving loads.

$$\delta t_{proper} \propto \frac{1}{V} \qquad (8)$$

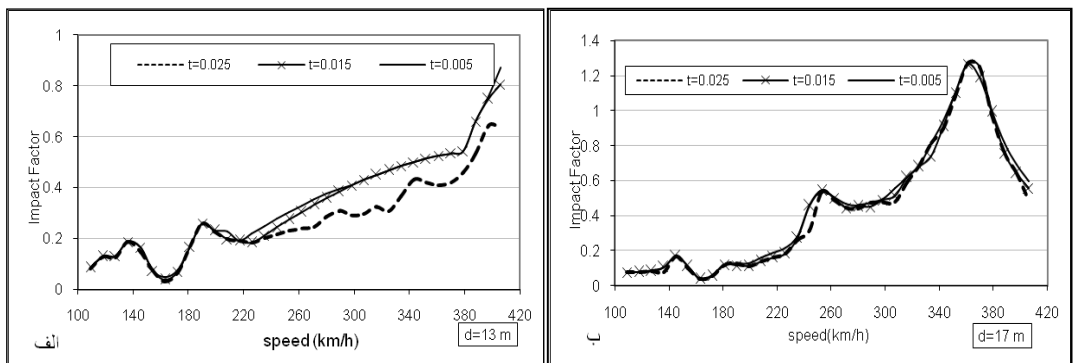

Figure 7: Impact factor for bridge with 10 m span length and different time steps

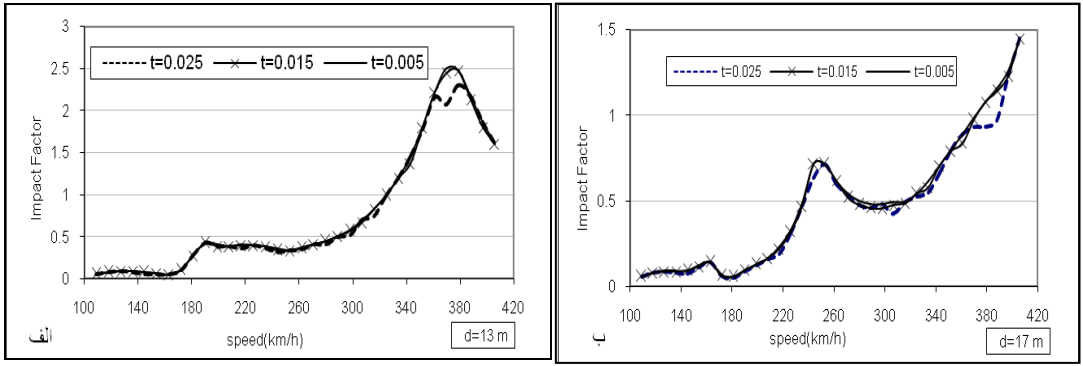
Figure 8: Impact factor for bridge with 15 m span length and different time steps

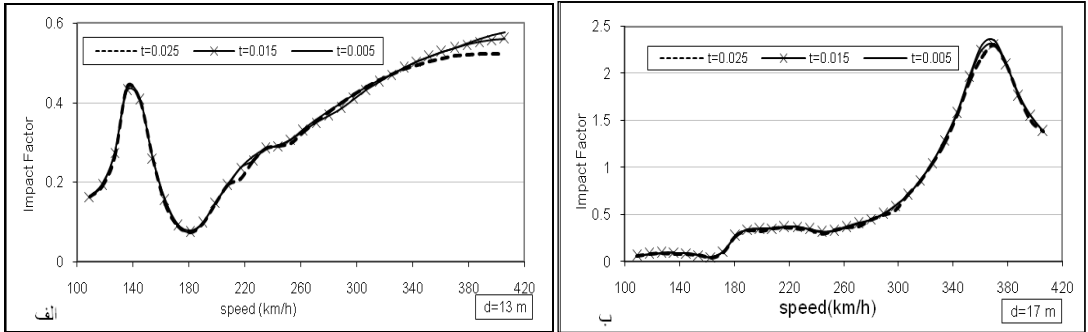
Figure 9: Impact factor for bridge with 20 m span length and different time steps

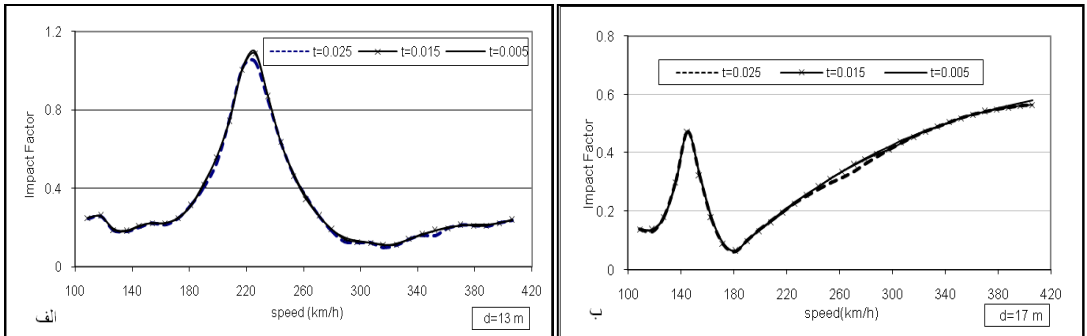
Figure 10: Impact factor for bridge with 25 m span length and different time steps

### *5-2- span length*

As shown in figures 7-10, the impact factor calculated by dynamic analyses with different time steps has smaller variation for bridges with longer span. This means according to dynamic analyses carried out in this research, it is possible to use biger time step for bridges with biger span length. In fact increase in span length causes increase in period of bridge viberation. Considering this fact in most of integration methods, time step is determined as a coefficient of period of structure, in this case too the proper time step for dynamic anaylsis of bridge under moving loads is proportional to span length (period of viberation).

$$\delta t_{proper} \propto L \qquad (9)$$

## 6- Determination of proper time step in dynamic analysis of bridges under moving loads

After determination of effective parameters on proper time step for dynamic analysis of bridges under moving loads, the equations 8 and 9 can be combined in equation 10:

$$\delta t_{proper} = k \frac{L}{V} \qquad (10)$$

In this equation L is span length in term of meter and V is velocity of load in meter per second and so the equation is dimensional compatible.

By comparison of proper time step for each condition with corresponding L/V the amount of k in equation 10 can be determined. For determination of proper time step from dynamic analysis it is necessary to specify an acceptable error for impact factor. As in different codes like AREMA and AASHTO, impact factor is calculated by two decimal [5, 10], so acceptable error for impact factor in this research is considered 0.01. This means if in a condition (span length and loading velocity) the difference of two impact factor calculated by two different time steps is less than 0.01 second, the bigger time step is considered as a proper time step for that condition. This method of determination of proper time step is shown in table 3 for a bridge with 10 meter span length and for velocity of motion. For example in velocity of 118 km/h, the difference of impact factors calculated by considering time step 0.05 , 0.025 second is more than 0.01 but this difference is less than 0.01 for time step 0.025 second and 0.015 second. So for this condition the proper time step is determined as a 0.025 second.

Table 3: Determination of proper time step according to dynamic analyses with different time steps

| speed(km/h) | Impact Factor | | | | | | Chosen time step |
|---|---|---|---|---|---|---|---|
| | δt=0.05 | δt=0.025 | Δt=0.015 | δt=0.01 | δt=0.005 | δt=0.0025 | |
| 109 | 0.08242 | 0.08242 | 0.08576 | 0.08576 | 0.08576 | 0.08596 | 0.05 |
| 118 | 0.10283 | 0.12949 | 0.12990 | 0.13222 | 0.13222 | 0.13222 | 0.025 |
| 145 | 0.14808 | 0.14808 | 0.16465 | 0.16010 | 0.16465 | 0.16465 | 0.015 |
| 208 | 0.19798 | 0.19798 | 0.19798 | 0.22182 | 0.22778 | 0.22919 | 0.01 |

After determination of time step for different conditions by method explained in table 3, and putting it in equation 10, the corresponding value for k is determined. In table 4 the minimum amount of k determined by this method is shown. The minimum amount of k for all 1632 conditions is 0.0195. Considering this amount for k in equation 10 in all condition is conservative. In table 4 in addition to minimum amount of k for some of conditions, minimum, Standard deviation and average of k for each bridge is shown. It is obvious the variation of k for different bridges is small and the result of this research is reliable for determination of proper time step is dynamic analysis of bridges under moving loads.

Table 4: Minimum values and variations of k

| axle distance (m) | Span length(m) | | | |
|---|---|---|---|---|
| | 10 | 15 | 20 | 25 |
| 13 | 0.0339 | 0.0351 | 0.0265 | 0.0205 |
| 15 | 0.0351 | 0.0209 | 0.0206 | 0.0302 |
| 17 | 0.0401 | 0.0251 | 0.0276 | 0.0281 |
| 22 | 0.0276 | 0.0527 | 0.0232 | 0.0328 |
| 23 | 0.0539 | 0.0284 | 0.0238 | 0.0291 |
| Minimum | 0.0276 | 0.0209 | 0.0206 | 0.0195 |
| Mean | 0.0362 | 0.0296 | 0.0240 | 0.0281 |
| Standard deviation | 0.0087 | 0.0113 | 0.0024 | 0.0061 |

As a final result the proper time step for dynamic analysis of bridges under moving loads is proposed by equation 11.

$$\delta t_{proper} = 0.195 \frac{L}{V} \qquad (11)$$

In this equation L is span length in term of meter and V is velocity of load in meter per second and the dimension of both sides is second.

## 7- Summery

In this research it is shown that although dynamic analysis of bridges under moving loads independent of time step size is stable and converges to results but the accuracy of results is dependent to the time step size. So it is important to determine the proper time step for dynamic analysis of bridges under moving loads. This requires implementation of a great number of dynamic analyses with different time steps which results to time and cost wasting. The equation proposed in this research for proper time step in dynamic analysis of bridges under moving loads prepare a basis for engineers and researchers to choose a reliable time step size in analyses without extra efforts, time and cost.

For investigation of effects of different parameter on proper time step for dynamic analysis of bridges under moving loads, dynamic responses of four bridges under different trains and velocities with different time steps were considered. When the difference of impact factors for a bridge for two time steps is smaller than 0.01 the bigger one chose as a proper time step. By this method for all 1603 conditions the proper time step was determined.

According to the results it is necessary to reduce the time step size by increase in load velocity. On the other hand by increase in span length which is equivalent to increase in flexural period of bridge it is possible to use bigger time step size in dynamic analysis. In this research by comparison between L/V of each condition and the chosen time step for that condition by dynamic analyses, a relation for proper time step in dynamic analysis of bridges under moving loads was proposed. This equation is the form of $\delta t_{proper} = k \dfrac{L}{V}$ and k is the minimum of all coefficients calculated in this research. This k is 0.0195 and conservative for all condition. For example if there was a proposed time step for dynamic analysis of bridge under moving loads, the number of dynamic analysis would be near one sixth of those carried out in this research.

## 8- References